\documentclass[graybox]{svmult}

\usepackage{amsmath}
\usepackage{type1cm}        
\usepackage{makeidx}         
\usepackage{graphicx}        
\usepackage{multicol}        
\usepackage[bottom]{footmisc}

\usepackage{newtxtext}       %
\usepackage{newtxmath}       
\usepackage{bbm}

\makeindex             

\begin{document}

\title*{The Information-Geometric Perspective of Compositional Data Analysis}
\author{Ionas Erb and Nihat Ay}
\authorrunning{Erb and Ay}
\titlerunning{The information-geometric perspective} 
\institute{Ionas Erb \at Centre for Genomic Regulation (CRG),
The Barcelona Institute of Science and Technology, Barcelona, Spain; \email{ionas.erb@crg.eu}
\and Nihat Ay \at Max-Planck Institute for Mathematics in the Sciences, Leipzig, Germany;\\
Department of Mathematics and Computer Science, Leipzig University, Leipzig, Germany;\\
Santa Fe Institute, Santa Fe, NM, USA; \email{nihat.ay@mis.mpg.de}}
%
%
\maketitle

\abstract{Information geometry uses the formal tools of differential geometry to describe the space of probability distributions as a Riemannian manifold with an additional dual structure. The formal equivalence of compositional data with discrete probability distributions makes it possible to apply the same description to the sample space of Compositional Data Analysis (CoDA). The latter has been formally described as a Euclidean space with an orthonormal basis featuring components that are suitable combinations of the original parts. In contrast to the Euclidean metric, the information-geometric description singles out the Fisher information metric as the only one keeping the manifold's geometric structure invariant under equivalent representations of the underlying random variables. Well-known concepts that are valid in Euclidean coordinates, e.g., the Pythagorean theorem, are generalized by information geometry to corresponding notions that hold for more general coordinates. In briefly reviewing Euclidean CoDA and, in more detail, the information-geometric approach, we show how the latter justifies the use of distance measures and divergences that so far have received little attention in CoDA as they do not fit the Euclidean geometry favoured by current thinking. We also show how Shannon entropy and relative entropy  can describe amalgamations in a simple way, while Aitchison distance requires the use of geometric means to obtain more succinct relationships. We proceed to prove the information monotonicity property for Aitchison distance. We close with some thoughts about new directions in CoDA where the rich structure that is provided by information geometry could be exploited.}

\section{Introduction}
Information geometry and Compositional Data Analysis (CoDA) are fields that have ignored each other so far. Independently, both have found powerful descriptions that led to a deeper understanding of the geometric relationships between their respective objects of interest: probability distributions and compositional data. Although both of these live on the same mathematical space (the simplex), and some of the mathematical structures are identically described, surprisingly, both fields have come to focus on quite different geometric aspects. On the one hand, the tools of differential geometry have revealed the underlying duality of the manifold of probability distributions, with the Fisher information metric playing a central role. On the other hand, the classical log-ratio approach led to the modern description of the compositional sample space as Euclidean and affine. We think it is  time that CoDA starts to profit from the rich structures that information geometry has to offer. This paper intends to build some bridges from information geometry to CoDA. In the first section, we will give a brief description of the Euclidean CoDA perspective. The second, and main, part of the paper describes in some detail the approach centred around the Fisher metric, with a description of the dual coordinates, exponential families and of how dual divergence functions generalize the notion of Euclidean distance. To ease understanding, throughout this section we link those concepts to the ones familiar in CoDA. In the third part of the paper, we show how information-based measures can lead to simpler expressions when amalgamations of parts are involved, and an important monotonicity result that holds for relative entropy is derived for Aitchison distance. We conclude with a short discussion about where we could go from here.      

\section{The Euclidean CoDA perspective}\label{CoDA}
Compositional data analysis is now unthinkable without the log-ratio approach pioneered by Aitchison \cite{AitchisonBook}. It has both led to a variety of data-analytic developments (for the most recent review, see \cite{green}), and more formal mathematical descriptions (see \cite{sampleSpace}). Following these, compositions can be described as equivalence classes whose representatives are points in a Euclidean space. We will give a brief recount here for completeness of exposition. Compositional data are defined as vectors of strictly positive numbers describing the parts of a whole for which the relevant information is only relative. As such, the absolute size of a $D$-part composition $\boldsymbol{x}\in \mathbb{R}^D_+$ is irrelevant, and all the information is conveyed by the ratios of its components. This can further be formalized by considering equivalent compositions $\boldsymbol{x}$, $\boldsymbol{y}$ such that $\boldsymbol{y}=c\boldsymbol{x}$ for a positive constant $c$. A composition is then an equivalence class of such proportional vectors \cite{mathCoDA}. A closed composition is the simplicial representative $\mathcal{C}\boldsymbol{x}:=\boldsymbol{x}/\sum_ix_i$, where the symbol $\mathcal{C}$ denotes the closure operation (i.e., the division by the sum over the components). Closed compositions are elements of the simplex
\begin{equation}
    \mathcal{S}^D:=\left\{(x_1,\dots,x_D)^T\in\mathbb{R}^D:x_i>0,i=1,\dots,D,\sum_i^Dx_i=1\right\},
\end{equation}
where $T$ denotes transposition. The simplex $\mathcal{S}^D$ can be equipped with a Euclidean structure by the vector space operations of perturbation and powering (playing the role of vector addition and scalar multiplication in real vector spaces), defined by
\begin{eqnarray}
\boldsymbol{x}\oplus\boldsymbol{y}&:=&\mathcal{C}(x_1y_1,\dots,x_Dy_D)^T,  \label{addition} \\
\alpha\odot\boldsymbol{x}&:=&\mathcal{C}(x_1^\alpha,\dots,x_D^\alpha)^T. \label{scalarmultiplication}
\end{eqnarray}
An inverse perturbation is given by $\ominus\boldsymbol{x}:=(-1)\odot\boldsymbol{x}$.\\
As a vector space, ${\mathcal S}^D$ also carries the structure of an affine space, and we can study affine subspaces, which are referred to as {\em linear manifolds\/} in \cite{simplicialGeometry}. In order to do so, we require a set of vectors $\boldsymbol{x}_1,\dots , \boldsymbol{x}_m$, which we assume to be linearly independent, and an origin point  $\boldsymbol{x}_0$. Here, independence means the following:  
Let $\boldsymbol{n}=\mathcal{C}(1,\dots,1)$ be the neutral element. A set of $m$ compositions $\boldsymbol{x}_1,\dots,\boldsymbol{x}_m\in\mathcal{S}^D$ is called {\it perturbation-independent} if the fact that $\boldsymbol{n}=\bigoplus_{i=1}^m(\alpha_i\odot\boldsymbol{x}_i)$ implies $\alpha_1=\dots=\alpha_m=0$. 
With this, an affine subspace is given as the set of compositions $\boldsymbol{y}$ such that
\begin{equation} \label{affsubspaces}
    \boldsymbol{y}=\boldsymbol{x}_0\oplus\bigoplus_{i=1}^m(\alpha_i\odot\boldsymbol{x}_i)
\end{equation}
for any real constants $\alpha_i$, $i=1,\dots,m$. Due to the linear independence of the vectors 
$\boldsymbol{x}_1,\dots , \boldsymbol{x}_m$, this is 
an $m$-dimensional space. \\
It is convenient to define the inner product for our Euclidean space via the so-called centred log-ratio transformation. Its definition \cite{logratioTrans} and inverse operation are 
\begin{eqnarray}
    \boldsymbol{v}&=&\mathrm{clr}(\boldsymbol{x}):=\left(\log\frac{x_1}{g(\boldsymbol{x})},\dots,\log\frac{x_D}{g(\boldsymbol{x})}\right)^T,\\
    \boldsymbol{x}&=&\mathrm{clr}^{-1}(\boldsymbol{v})=\mathcal{C}\mathrm{exp}(\boldsymbol{v}),
\end{eqnarray}
where $g$ denotes the geometric mean 
$g(\boldsymbol{x}) = \left(\prod_{i = 1}^D x_i\right)^{\frac{1}{D}}$. Note that the sum over the components $\mathrm{clr}_i$ of clr-transformed vectors is 0.  
An inner product can then be defined by
\begin{equation}
    \left<\boldsymbol{x},\boldsymbol{y}\right>_A:=\sum_{i=1}^D\mathrm{clr}_i(\boldsymbol{x})\mathrm{clr}_i(\boldsymbol{y}).
\end{equation}
The corresponding (squared) norm and distance are
\begin{equation}
    \left\lVert\boldsymbol{x}\right\rVert^2_A=\left<\boldsymbol{x},\boldsymbol{x}\right>_A,~~~~~d^2_A(\boldsymbol{x},\boldsymbol{y})=\left\lVert\boldsymbol{x}\ominus\boldsymbol{y}\right\rVert^2_A.
\end{equation}
This distance is known as Aitchison distance (denoted by the $A$ subscript). Note that the clr-transformation is an isometry $\mathcal{S}^D\rightarrow \mathcal{T}^D \subset\mathbb{R}^D$ between 
$(D-1)$-dimensional Euclidean spaces, 
where 
\begin{equation}
\mathcal{T}^D :=\left\{\boldsymbol{v}\in\mathbb{R}^D:\sum_{i=1}^Dv_i=0\right\}.
\end{equation}
We can thus obtain orthonormal bases in $\mathcal{S}^D$ from orthonormal bases in $\mathcal{T}^D$. Such orthonormal basis vectors are defined by the columns 
$\boldsymbol{v}_i$, $i = 1,\dots, D-1$, of $\mathbf{V}$, a matrix of order $D\times (D-1)$ obeying
\begin{eqnarray}
\mathbf{V}^T\mathbf{V}&=&\mathbf{I}_{D-1},\label{contrast1}\\
\mathbf{V}\mathbf{V}^T&=&\mathbf{I}_D-\frac{1}{D}\mathbf{1}\mathbf{1}^T,\label{contrast2}
\end{eqnarray}
with $\mathbf{I}_D$ the $D\times D$ identity matrix and $\mathbf{1}$ a $D$-dimensional vector where each component is 1. The first equation ensures orthonormality, the second equation makes sure components sum to zero. Now the vectors
\begin{equation}
    \boldsymbol{e}_i=\mathrm{clr}^{-1}(\boldsymbol{v}_i)=\mathcal{C}\mathrm{exp}(\boldsymbol{v}_i)\label{basis}
\end{equation}
constitute an orthonormal basis in $\mathcal{S}^D$. Euclidean coordinates $\boldsymbol{z}$ of $\boldsymbol{x} \in \mathcal{S}^D$ with respect to the basis 
$\boldsymbol{e}_i$, $i = 1,\dots, D-1$, then follow from the so-called isometric log-ratio transformation \cite{ilr}. Its definition and inverse operation are
\begin{eqnarray}
    \boldsymbol{z}=\mathrm{ilr}(\boldsymbol{x})&:=&\mathbf{V}^T\log(\boldsymbol{x}),\label{ilr}\\
    \boldsymbol{x}=\mathrm{ilr}^{-1}(\boldsymbol{z})&=&\mathcal{C}\mathrm{exp}(\mathbf{V}\boldsymbol{z}).
\end{eqnarray}
The second equation shows the composition as generated from the coordinates $\boldsymbol{z}$. This can also be written in the usual component form using the basis vectors of Eq.\ (\ref{basis}):
 \begin{equation}
     \boldsymbol{x}=\bigoplus_{i=1}^{D-1}(z_i\odot\boldsymbol{e}_i).
 \end{equation}
 Equations (\ref{contrast1}) and (\ref{contrast2}) characterize any orthonormal basis of $\mathcal{S}^D$. $\mathbf{V}$ is known under the name of contrast matrix. Many choices for this matrix are possible. Balance coordinates \cite{balances} are popular for their relative simplicity and interpretability, where pivot coordinates \cite{pivot} are sometimes preferred.

\section {The information-geometric perspective}

Information geometry started out as a study of the geometry of statistical estimation. The set of probability distributions that constitute a statistical model is seen as a manifold whose invariant geometric structures \cite{Chentsov} are studied in relation to the statistical estimation using this model. A Riemannian metric and a family of affine connections are naturally introduced on such a manifold \cite{Rao}. It turns out that a fundamental duality underlies these structures \cite{AmariNagaoka}. While the first ideas about the geometry of statistical estimation are from the first half of the 20th century and go back to a variety of authors, a first attempt at a unified exposition of the topic by Amari \cite{lectureNotes} was published around a similar time as Aitchison's book.  

Here we try to highlight some of the main concepts, borrowing from chapter two in \cite{Nihat}, but mainly following the treatment in \cite{Amari}. The latter also showcases the many applications information geometry has found during the last decades. While the best known book on the topic may be \cite{AmariNagaoka}, the most formal and complete exposition can currently be found in \cite{Nihat}.\\
The ideas that require more advanced notions from differential geometry will not be touched upon in our short outline.

\subsection{Dual coordinates and Fisher metric in finite information geometry }

To exploit the equivalence of probability distributions with compositional data, we only consider the manifold of discrete distributions $\mathcal{S}^D$. To emphasize the equivalence, we will denote the probabilities by the same symbol as our compositions, i.e., $\boldsymbol{x}\in\mathcal{S}^D$ is a vector of $D$ probabilities. 
To complete the probabilistic picture, we need a random variable $R$ which can take values $r \in \{1,\dots,D\}$ with the respective probabilities $x_r = \mathbb{P}\{R = r\}$, the coordinates of $\boldsymbol{x}$.
The distribution of this random variable can now be written as
\begin{equation} \label{repdist}
    p_R = \boldsymbol{x} = (x_1,\dots, x_D)^T \in \mathcal{S}^D.
\end{equation}

From the information-geometric perspective, there are two natural ways to parametrize the set ${\mathcal S}^D$ of all (strictly positive) distributions of $R$. The first possibility is quite obvious. 
The first $D-1$ probabilities in Eq.\ (\ref{repdist}) that are free to be specified, that is $x_1,\dots, x_{D-1}$, 
can be considered parameters, which we denote by $\boldsymbol{\eta}=(x_1,\dots,x_{D-1})^T$. 
In these coordinates, probability distributions are written as
\begin{equation} \label{etapara}
    p(r ; {\boldsymbol{\eta}}) =
    \left\{
      \begin{array}{c@{,\quad}l}
         \eta_r & \mbox{if $r \leq D-1 $} \\
         1 - \sum_{i = 1}^{D-1} \eta_i 
         & \mbox{if $r = D$}
      \end{array}
    \right. , \qquad r = 1,\dots,D.
\end{equation}

Alternatively, our distribution can be parametrized using what is known as the alr-transformation \cite{logratioTrans} in CoDA:
\begin{equation}
    \theta^i=\log\frac{x_i}{x_D},~~~i=1,\dots,D-1.
\end{equation}
With this, we can write our distribution in the form  
\begin{equation} 
    p(r; {\boldsymbol{\theta}}) =  \mathrm{exp}\left(\sum_{k=1}^{D-1}\theta^k {\mathbbm 1}_k(r) -\psi(\boldsymbol{\theta})\right), \qquad 
    r = 1,\dots,D,
    \label{thetadist2}
\end{equation}
where ${\mathbbm 1}_k(r) = 1$ if $r = k$, and ${\mathbbm 1}_k(r) = 0$ otherwise.   
The function $\psi$ ensures normalization, that is
\begin{equation} \label{psi}
    \psi(\boldsymbol{\theta})=\log\left(1+\sum_{i=1}^{D-1}e^{\theta^i}\right)=-\log x_D.
\end{equation}
{The parametrization of Eq.\ (\ref{thetadist2}) in terms of $\boldsymbol{\theta} \in \mathbb{R}^{D-1}$ can be used in order to define a linear structure on ${\mathcal S}^D$. The addition of two 
distributions $p(\cdot ; \boldsymbol{\theta}_1)$ and $p(\cdot ; \boldsymbol{\theta}_2)$ can simply be defined by taking their product and then normalizing. With this vector addition, denoted by $\oplus$, we obviously 
have     
\begin{eqnarray}
    p({r}; \boldsymbol{\theta}_1) \oplus p({r}; \boldsymbol{\theta}_2) 
    & = & p({r}; \boldsymbol{\theta}_1 + \boldsymbol{\theta}_2).
\end{eqnarray}
Thus, the operation $\oplus$ is consistent with the usual addition in the parameter space $\mathbb{R}^{D - 1}$. 
The multiplication of a distribution $p(\cdot ; \boldsymbol{\theta})$ with a scalar $\alpha \in \mathbb{R}$ can be correspondingly defined by potentiating and then normalizing. This defines a scalar multiplication $\odot$, and we have      
\begin{eqnarray} 
   {\alpha} \odot p({r} ; \boldsymbol{\theta})
    & = & p({r} ; \alpha \cdot \boldsymbol{\theta}).
\end{eqnarray}
Obviously, the scalar multiplication $\odot$ is consistent with the usual multiplication in the parameter space. Note that the vector space structure defined here, which is well known in information geometry, coincides with the structure defined by equations (\ref{addition}) and (\ref{scalarmultiplication}). Given the linear structure, we can consider affine subspaces of $\mathcal{S}^D$. These are well-known and fundamental families in statistics, statistical physics, and information geometry, the so-called {\em exponential families\/}. We basically obtain them from the representation of Eq.\ (\ref{thetadist2}) if we replace the functions ${\mathbbm 1}_k$ by $d$ ($d\le D-1$) functions $X_k: \{1,\dots,D\} \to \mathbb{R}$, and shift the whole family by some reference measure $p_0$: 
\begin{equation} 
    p(r; {\boldsymbol{\theta}}) =  p_0(r) \,  \mathrm{exp}\left(\sum_{k=1}^{d}\theta^k X_k(r) -\psi(\boldsymbol{\theta})\right), \qquad 
    r = 1,\dots,D.
    \label{thetadist3}
\end{equation}
Here, $\psi(\boldsymbol{\theta})$ ensures normalization,  
but it does not reduce to the simple structure of Eq.\ (\ref{psi}) in general.  
In statistical physics, the function $\psi$ 
is known as the free energy (in other contexts it is also known as the cumulant-generating function).
Note that, given the same linear structure on $\mathcal{S}^D$, the 
exponential families coincide with the linear manifolds of Eq.\ (\ref{affsubspaces}), which were introduced into the field of Compositional Data Analysis more recently.}  \\
\medskip

\begin{figure}[t]
\includegraphics[width=12cm]{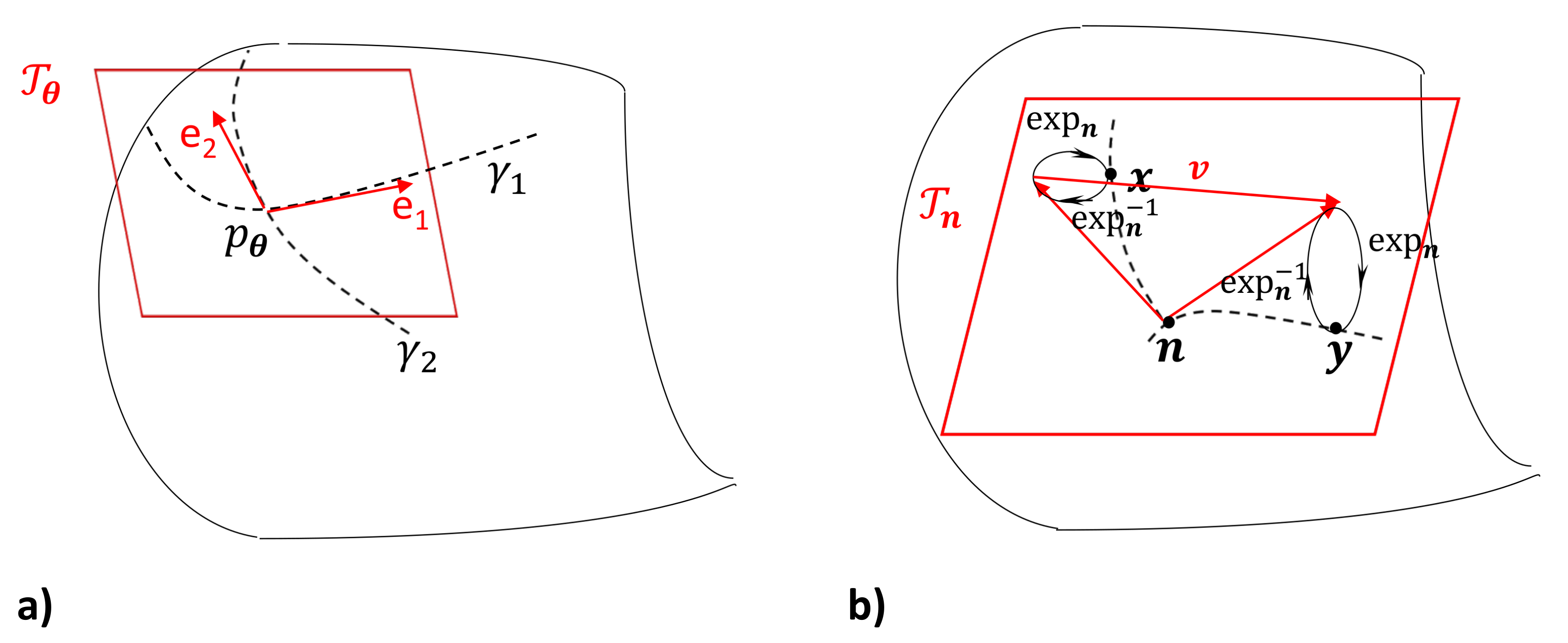}
\caption{{\bf a)} A manifold (in black) and its tangent space in point $p_\theta$ (in red) where two coordinate curves (dashed lines) are crossing. The velocities of the coordinate curves are the basis vectors $\mathrm{e_1}$ and $\mathrm{e_2}$. {\bf b)} Exponential map and its inverse via the tangent space anchored in a reference point $\boldsymbol{n}$. The points $\boldsymbol{x}$ and $\boldsymbol{y}$ are mapped via $\mathrm{exp}_{\boldsymbol{n}}^{-1}$ to the tangent space. The difference vector between the two mapped points is $\boldsymbol{v}=\mathrm{vec}(\boldsymbol{x},\boldsymbol{y})$.}
\label{manifold}
\end{figure}

In what follows we restrict attention to the parametrizations of equations (\ref{etapara}) and (\ref{thetadist2}) of the full simplex $\mathcal{S}^D$ as one instance of the general structure that underlies information geometry. The function $\psi$, given by Eq.\ (\ref{psi}), is a convex function, and we can get back the $\boldsymbol{\eta}$ coordinates from it via a Legendre transformation \cite{AmariNagaoka}
\begin{equation}
    \boldsymbol{\eta}=\nabla\psi(\boldsymbol{\theta}),\label{eta}
\end{equation}
where $\nabla$ denotes the partial derivatives $(\partial\psi/\partial\theta^i)_{i=1}^D$.
The Legendre dual of $\psi(\boldsymbol{\theta})$ is another convex function defined by
\begin{equation}
    \phi(\boldsymbol{\eta})=\max_{\boldsymbol{\theta}}\left\{\boldsymbol{\theta}\cdot\boldsymbol{\eta}-\psi(\boldsymbol{\theta})\right\},
\end{equation}
which is given by the negative Shannon entropy
\begin{equation}
    \phi(\boldsymbol{\eta})=\sum_{i=1}^{D-1}\eta_i\log\eta_i+\left(1-\sum_{i=1}^{D-1}\eta_i\right)\log\left(1-\sum_{i=1}^{D-1}\eta_i\right).
\end{equation}
 In equivalence to Eq.\ (\ref{eta}), from $\phi(\boldsymbol{\eta})$ we can get back the $\boldsymbol{\theta}$ coordinates:
 \begin{equation}
    \boldsymbol{\theta}=\nabla\phi(\boldsymbol{\eta}).
\end{equation}
 There is thus a fundamental duality mediated by the Legendre transformation which links the two types of parameters $\boldsymbol{\eta}$ and $\boldsymbol{\theta}$ as well as the convex functions $\psi$ and $\phi$. Legendre transformations are well known to play a fundamental role in phenomenological thermodynamics. Their importance for information geometry was established by Amari and Nagaoka \cite{AmariNagaoka}.\\
In what follows, we use the parameters $\boldsymbol{\theta}$. Like $\boldsymbol{\eta}$, they define a (local) coordinate system of the manifold $\mathcal{S}^D$. From each point $\boldsymbol{\theta}$, $D-1$ coordinate curves $\gamma_i: t \mapsto \gamma_i(t)$ in ${\mathcal S}^D$,
$\gamma_i(0) = p_{\boldsymbol{\theta}}$, 
emerge when holding $D-2$ of the $\theta_i$ constant. 
Consider their velocities  
\begin{equation} 
    \mathrm{e}_i := 
    \left. \frac{d}{dt}\gamma_i(t) \right|_{t = 0} , \qquad i = 1,\dots, D-1. \label{tangent}
\end{equation}
In the point $\boldsymbol{\theta}$ itself, the $D-1$ vectors $\mathrm{e}_i$ pointing in the direction of each coordinate curve form the basis of the so-called {\it tangent space} ${T}_{\boldsymbol{\theta}} \mathcal{S}^D = {\mathcal T}^D$ of this point, see Figure \ref{manifold}a). Similarly, we can define vectors $\mathrm{e}^i$ with respect to the coordinates $\boldsymbol{\eta}$, which also span the tangent space.
\medskip

For a Riemannian manifold a metric tensor $\mathbf{G}$ is defined. This metric is obtained via an inner product on the tangent space:
\begin{equation}
    g_{ij}=\left<\mathrm{e}_i,\mathrm{e}_j\right>,
\end{equation}
which depends on the coordinates chosen. The coordinate system is Euclidean if $g_{ij}=\delta_{ij}$ (this is the case for the parametrization achieved by Eq.\ \ref{ilr}). 
For the Riemannian metric of exponential families, the basis vectors can be identified with the so-called score\footnote{The score function plays an important role in maximum-likelihood estimation.} function $\partial\log p(r,\boldsymbol{\theta})/\partial\theta^i$. The resulting Riemannian metric is known as the Fisher information matrix:
\begin{eqnarray}
    g_{ij}(\boldsymbol{\theta})
    &=& E_{\boldsymbol{\theta}} \left[\frac{\partial}{\partial\theta^i}\log p({r};\boldsymbol{\theta})\, \frac{\partial}{\partial\theta^j}\log p({r};\boldsymbol{\theta}) \right], \\
    & = & E_{\boldsymbol{\theta}}\left[ (\mathbbm{1}_i - E(\mathbbm{1}_i))(\mathbbm{1}_j - E(\mathbbm{1}_j))\right] \label{cova} \\
    & = & 
    \left\{ 
       \begin{array}{c@{,\quad}l}
          - \eta_i \eta_j & \mbox{if $i \not= j$} \\
          \eta_i (1 - \eta_i) & \mbox{if $i = j$}
       \end{array}
    \right. .
\end{eqnarray}
with $E$ denoting the expectation value, and  
\[
  \eta_i = \frac{e^{\theta^i}}{1 + \sum_{k = 1}^{D-1} e^{\theta^k}}, \qquad i = 1, \dots, D-1.
\]
Note that this is the covariance matrix of the random vector 
$(\mathbbm{1}_1, \dots, \mathbbm{1}_{D - 1})$, as expressed by Eq.\ (\ref{cova}). 
Convex functions have positive definite Hessian matrices. Here their elements are given by the Fisher metric itself:
\begin{eqnarray}
    g_{ij}(\boldsymbol{\theta})&=&\frac{\partial}{\partial\theta^i}\frac{\partial}{\partial\theta^j}\psi(\boldsymbol{\theta}),\\
    g^{ij}(\boldsymbol{\eta})&=&\frac{\partial}{\partial\eta_i}\frac{\partial}{\partial\eta_j}\phi(\boldsymbol{\eta}).
\end{eqnarray}
The second matrix is the inverse of the first, which follows from their Legendre duality. This means the second matrix is an inverse covariance, an important object in the theory of graphical models \cite{graphMod}. Although the Fisher metric is not Euclidean, i.e., $\left<\mathrm{e}_i,\mathrm{e}_j\right>\ne\delta_{ij}$, we do have a generalization of this when mixing the dual coordinates: $\left<\mathrm{e}_i,\mathrm{e}^j\right>=\delta_{ij}$. \\

Note that for a general exponential family of the form of Eq.\ (\ref{thetadist3}), the Fisher information matrix can be expressed as a covariance matrix, which generalizes Eq.\ (\ref{cova}):  
\begin{eqnarray}
    g_{ij}(\boldsymbol{\theta})
    &=& E_{\boldsymbol{\theta}} \left[\frac{\partial}{\partial\theta^i}\log p({r};\boldsymbol{\theta})\, \frac{\partial}{\partial\theta^j}\log p({r};\boldsymbol{\theta}) \right], \\
    & = & E_{\boldsymbol{\theta}}\left[ (X_i - E(X_i))(X_j - E(X_j))\right]. \label{cova2}
\end{eqnarray}

\subsection{Distance measures and divergences}\label{divergences}

\begin{figure}[t]
\sidecaption[t]
\includegraphics[width=7.5cm]{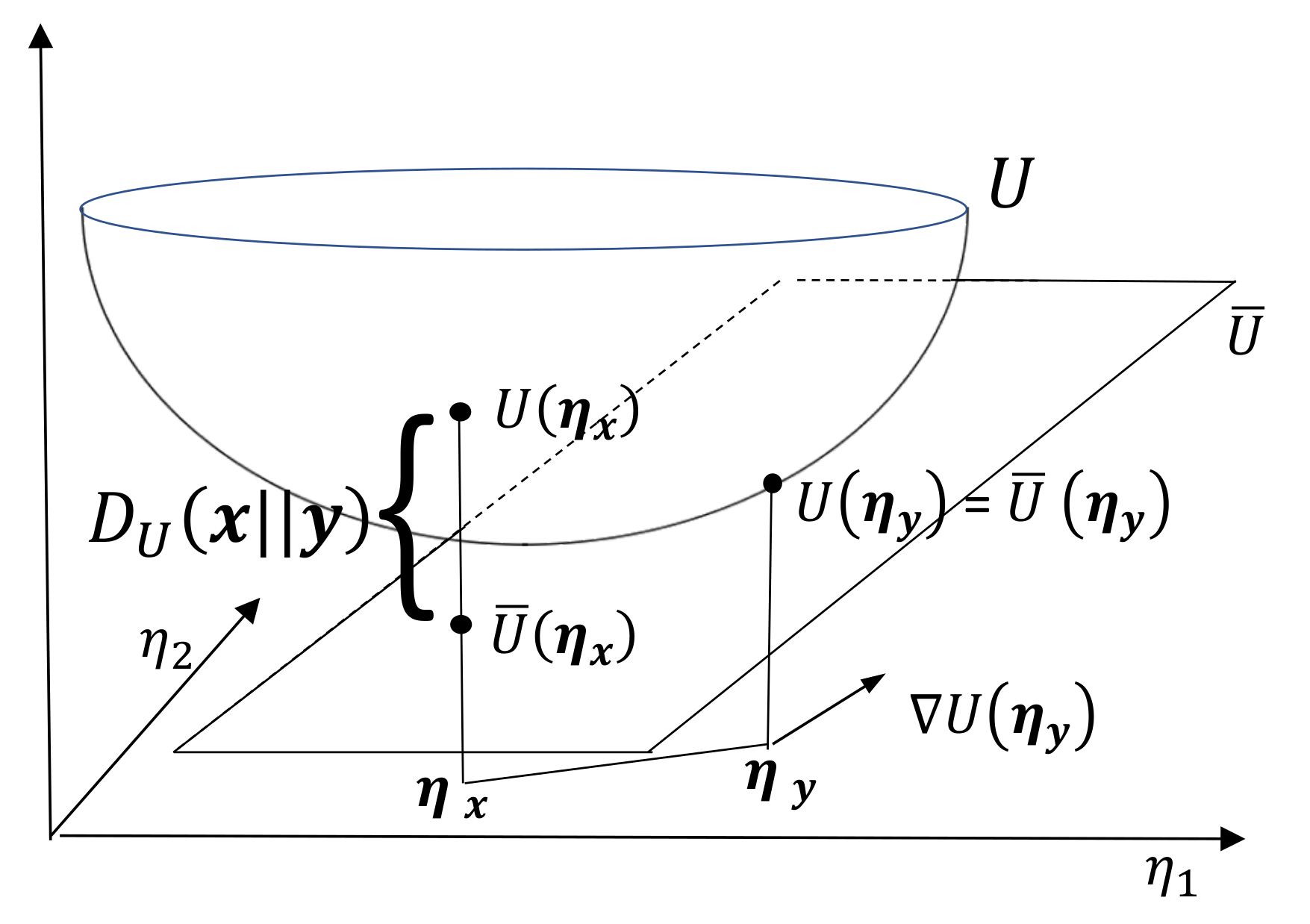}
\caption{The graph of the potential $U$, of its linearization $\overline{U}$, and a visualization of the corresponding Bregman divergence $D_U$. The gradient vector $\nabla U$ points in the direction of greatest change of the function $U$ in the space of the coordinates $\boldsymbol{\eta}$.} 
\label{bregman}
\end{figure}

Our affine structure on $\mathcal{S}^D$ can be reformulated additively via an exponential map $\mathcal{S}^D\times\mathcal{T}^D\to\mathcal{S}^D $, $(\boldsymbol{x},\boldsymbol{v})\mapsto\boldsymbol{y}$: 
\begin{eqnarray}
    \boldsymbol{y}=\mathrm{exp}_{\boldsymbol{x}}(\boldsymbol{v})&:=&\boldsymbol{x}\oplus e^{\boldsymbol{v}}, \label{aff1}\\
    \boldsymbol{v}=\mathrm{vec}(\boldsymbol{x},\boldsymbol{y})&:=&\mathrm{exp}^{-1}_{\boldsymbol{x}}(\boldsymbol{y})=\mathrm{clr}(\boldsymbol{y})-\mathrm{clr}(\boldsymbol{x}),\label{aff2}
\end{eqnarray}
using the notation introduced in section \ref{CoDA}, and shown here together with its inverse. This map is used in \cite{AyErb}, where ordinary linear differential equations with the time derivative defined for $\mathrm{vec}(\boldsymbol{x}(t_0),\boldsymbol{x}(t))$ letting $t\to t_0$ are considered. These turn out to be replicator equations with special properties that are known from population dynamics. Note that, e.g., with the center of the simplex $\boldsymbol{n}$ as defined before, we have $\mathrm{vec}(\boldsymbol{n},\boldsymbol{x})=\mathrm{clr}(\boldsymbol{x})$. With the exponential map, we can interpret $\mathrm{vec}(\boldsymbol{x},\boldsymbol{y})$ as the difference vector between two compositions, see Figure \ref{manifold}b). The set of all such difference vectors for a given point can be interpreted as the gradient field of a convex function 
(a.k.a.\ a potential). 
In order to highlight the generality of this concept in information geometry, we consider a general convex function $U$ with respect to some parameters $\boldsymbol{\eta}$ from a convex domain. In the following, we will denote the compositions for which the parameters are to be evaluated by subscripts. The linearization of $U$ in $\boldsymbol{\eta}_{\boldsymbol{y}}$ is 
given by 
\[
    \overline{U}(\boldsymbol{\eta}_{\boldsymbol{x}}) \, := \, U(\boldsymbol{\eta}_{\boldsymbol{y}}) + 
            \nabla U (\boldsymbol{\eta}_{\boldsymbol{y}}) \cdot (\boldsymbol{\eta}_{\boldsymbol{x}} - \boldsymbol{\eta}_{\boldsymbol{y}}).
\]
The graph of this linearization is a hyperplane of dimension $D-1$ touching the graph of $U$ in the point 
$(\boldsymbol{\eta}_{\boldsymbol{y}}, U(\boldsymbol{\eta}_{\boldsymbol{y}}))$, see Figure \ref{bregman}. The difference between $U$ and its linearisation $\overline{U}$ in $\boldsymbol{\eta}_{\boldsymbol{y}}$, evaluated at $\boldsymbol{\eta}_{\boldsymbol{x}}$, defines a so-called {\em Bregman divergence\/}, a class of divergences that plays an important role in information geometry. More precisely,  
\begin{equation}
    D_U(\boldsymbol{x}||\boldsymbol{y}) \, := \,  
    U(\boldsymbol{\eta}_{\boldsymbol{x}} ) - \overline{U}(\boldsymbol{\eta}_{\boldsymbol{x}})
    \, = \, U(\boldsymbol{\eta}_{\boldsymbol{x}})-U(\boldsymbol{\eta}_{\boldsymbol{y}})-\nabla U(\boldsymbol{\eta}_{\boldsymbol{y}})\cdot(\boldsymbol{\eta}_{\boldsymbol{x}}-\boldsymbol{\eta}_{\boldsymbol{y}}).
    \label{Bregman}
\end{equation}

Divergences are similar to distance functions but they are not necessarily symmetric and need not fulfill the triangle inequality. As an example, 
let us consider the potential naturally associated with the structure of equations (\ref{aff1}) and (\ref{aff2}), the squared Aitchison norm
\begin{equation}
    A(\boldsymbol{z}_{\boldsymbol{x}})=\sum_{i=1}^D\mathrm{clr}_i(\boldsymbol{x})^2=\sum_{i=1}^{D-1}\mathrm{ilr}_i(\boldsymbol{x})^2,\label{ilrpot}
\end{equation}
with ilr$_i$ the $i$-th ilr coordinate $z_i$, see Eq.\ (\ref{ilr}), and $\boldsymbol{z}_{\boldsymbol{x}}$ the vector of coordinates $z_i$. We then have
\begin{equation}
    D_A(\boldsymbol{x}||\boldsymbol{y})=\sum_{i=1}^{D-1}(\mathrm{ilr}_i(\boldsymbol{x})-\mathrm{ilr}_i(\boldsymbol{y}))^2=d^2_A(\boldsymbol{x},\boldsymbol{y}),\label{EuclideanDiv}
\end{equation}
coinciding with the squared Aitchison distance. This is the special case of a Euclidean divergence, which is also a (squared) distance function. 

Let us now come to the divergences that arise when replacing $U(\boldsymbol{\eta})$ in Eq.\ (\ref{Bregman}) by our dually convex functions $\psi(\boldsymbol{\theta})$ and $\phi(\boldsymbol{\eta})$. They turn out to be the relative entropies (a.k.a.\ Kullback-Leibler divergences)
\begin{eqnarray}
    D_\phi(\boldsymbol{x}||\boldsymbol{y})&=&\sum_{i=1}^Dx_i\log\frac{x_i}{y_i},\\
    D_\psi(\boldsymbol{x}||\boldsymbol{y})&=&\sum_{i=1}^Dy_i\log\frac{y_i}{x_i}.
\end{eqnarray}
Thus the symmetry we had in the Euclidean case finds its generalization for our dual case in 
\begin{equation}
  D_\psi(\boldsymbol{x}||\boldsymbol{y})=D_\phi(\boldsymbol{y}||\boldsymbol{x}).  
\end{equation}
Moreover, one can show that we can ``complete the square'' via
\begin{equation}
    D_\psi(\boldsymbol{x}||\boldsymbol{y})=\psi(\boldsymbol{\theta}_{\boldsymbol{x}})+\phi(\boldsymbol{\eta}_{\boldsymbol{y}})-\boldsymbol{\theta}_{\boldsymbol{x}}\cdot\boldsymbol{\eta}_{\boldsymbol{y}}.
\end{equation}
Of course, symmetrizations of relative entropy exist, with the Jensen-Shannon divergence perhaps the most prominent example. Also, a symmetric ``compositional'' version of relative entropy has been proposed because of its additional properties that are often regarded as indispensable\footnote{As an example, the translation invariance of distance measures, known under the name of perturbation invariance in CoDA, has its information-geometric analog in the invariance of the inner product of tangent vectors $\mathbf{u}$, $\mathbf{v}$ under parallel transport $P$ and its dual $P^*$: $\left<\mathbf{u},\mathbf{v}\right>=\left<P(\mathbf{u}),P^*(\mathbf{v})\right>$.} in CoDA \cite{symmDiv}. While such measures have some interesting properties, they do not make use of the duality of our parametrizations and are therefore less suitable for our approach. Indeed, although a dual divergence is not a distance measure in the strict sense, it can quantify the distance between points along a curve in a similar way, and generalizations of well-known relationships from Euclidean geometry are available. Geodesic lines connecting two compositions can be constructed via convex combinations of the parameters $\boldsymbol{\theta}$, and the corresponding dual geodesics from convex combinations of the dual parameters $\boldsymbol{\eta}$. Such geodesics are orthogonal to each other when the inner product 
of their tangent vectors, with respect to the Fisher metric,  vanishes in the point of intersection. In this case, a generalized Pythagorean theorem holds for the corresponding dual divergence, e.g.,
\begin{equation}
    D_\psi(\boldsymbol{x}||\boldsymbol{y})=D_\psi(\boldsymbol{x}||\boldsymbol{z})+D_\psi(\boldsymbol{z}||\boldsymbol{y}),
\end{equation}
where the geodesics $t\boldsymbol{\eta}_{\boldsymbol{x}}+(1-t)\boldsymbol{\eta}_{\boldsymbol{z}}$ and $t\boldsymbol{\theta}_{\boldsymbol{z}}+(1-t)\boldsymbol{\theta}_{\boldsymbol{y}}$ intersect orthogonally.

\subsection{Distances obtained from the Fisher metric compared with those used in CoDA}
\begin{figure}[t]
\sidecaption[t]
\includegraphics[width=6.5cm]{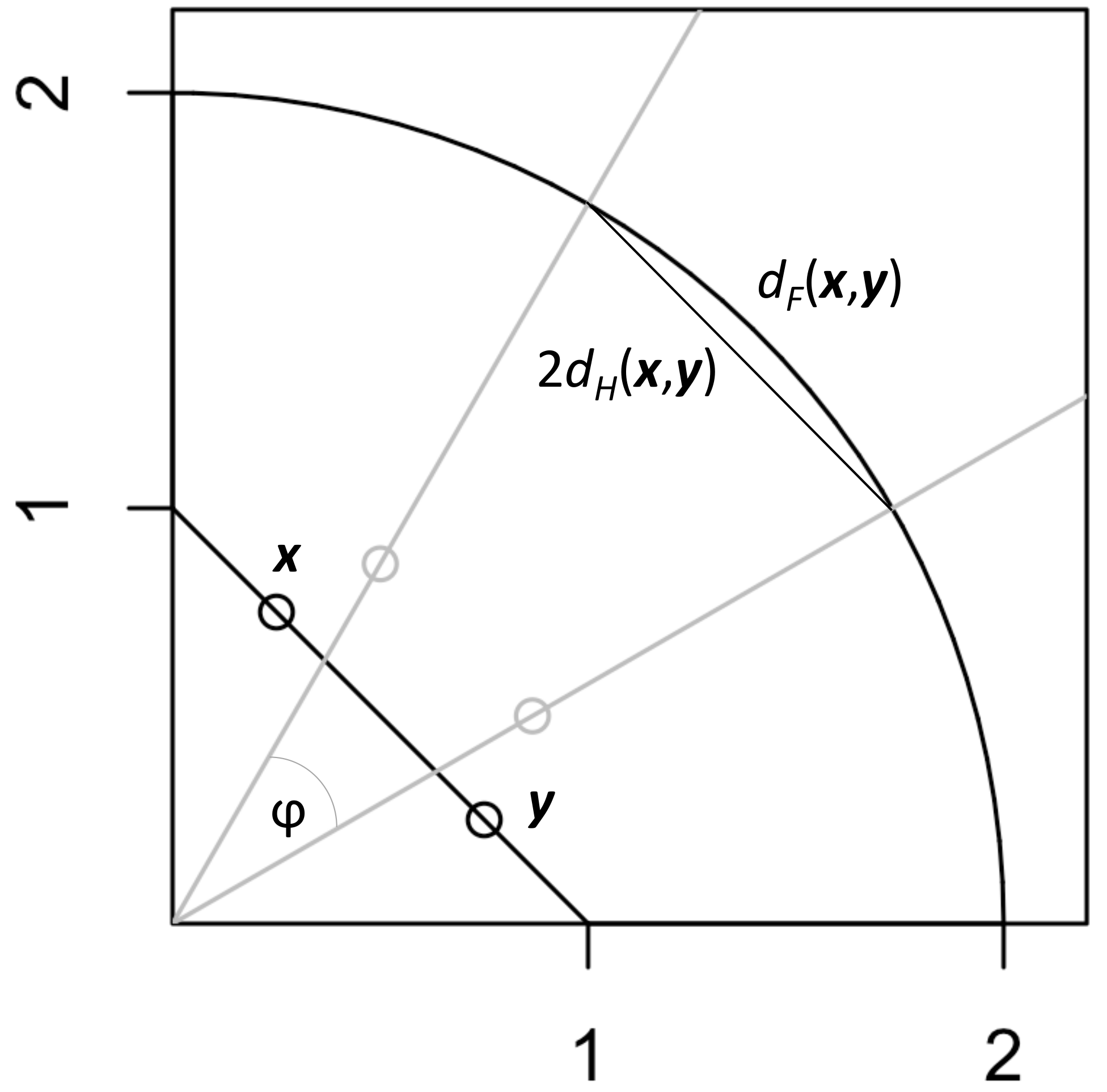}
\caption{Hellinger and Fisher distance between compositions $\boldsymbol{x}$ and $\boldsymbol{y}$. Shown are the 2-simplex and the positive orthant of a sphere with radius 2. The gray points indicate the square-root-transformed compositions. They lie on a sphere with radius one, and their Euclidean distance is the Hellinger distance. The lines from the origin going through these points cross the sphere with radius two. The length of the arc between the resulting points is the Fisher distance, while their Euclidean distance is twice the Hellinger distance.}
\label{hellinger}
\end{figure}

We have seen that the potential of Eq.\ (\ref{ilrpot}) led to a divergence that is also a Euclidean distance. Let us now consider a generalization of our dual divergences that includes as a special case a Euclidean distance that is related to the Fisher metric. The so-called $\alpha$-divergence (closely related to the Renyi \cite{Renyi} and Tsallis \cite{Tsallis} entropies) is defined as
\begin{equation}
    D_\alpha(\boldsymbol{x}||\boldsymbol{y})=\frac{4}{1-\alpha^2}\left(1-\sum_{i=1}^Dy_i^\frac{1+\alpha}{2}x_i^\frac{1-\alpha}{2}\right).\label{alphadiv}
\end{equation}
In the limit, the cases $\alpha=\pm1$ correspond to $D_\psi$ and $D_\phi$. The case $\alpha=0$ (where the divergence is self-dual, i.e., symmetric) corresponds to the Euclidean distance of the points $(\sqrt{x_1}, \dots, \sqrt{x_D})$ and $(\sqrt{y_1}, \dots, \sqrt{y_D})$:
\begin{eqnarray}
    d^2_H(\boldsymbol{x},\boldsymbol{y})\nonumber
    & = & \sum_{i=1}^D(\sqrt{x_i}-\sqrt{y_i})^2 \nonumber\\
    & = & \sum_{i=1}^D(x_i - 2 \sqrt{x_i y_i} + y_i) \\
    & = & 2 \left( 1 - 
        \sum_{i=1}^D \sqrt{x_i y_i} \right) \nonumber\\
    & = & \frac{1}{2} \, D_0 (\boldsymbol{x} \| \boldsymbol{y}).    
\end{eqnarray}
    This is the so-called Hellinger distance. It is closely related to the Riemannian distance\footnote{The Riemannian distance between two points on a manifold is the minimum of the lengths of all the piecewise smooth paths joining the two points.} between two compositions with respect to the Fisher metric. This so-called Fisher distance can be expressed explicitly \cite{Nihat} as
\begin{equation}
    d^2_F(\boldsymbol{x},\boldsymbol{y})=
    4 \, \mathrm{arccos}^2\left(\sum_{i=1}^D\sqrt{x_iy_i}\right),\label{Fisher}
\end{equation}
and its relation to the Hellinger distance is given by
\begin{equation}
    d^2_H(\boldsymbol{x},\boldsymbol{y})=2\left(1-\mathrm{cos}\frac{d_F(\boldsymbol{x},\boldsymbol{y})}{2}\right).
\end{equation}
These distances can be better understood when noting the role that plays the angle between the two points on the sphere, i.e., when comparing Eq.\ (\ref{Fisher}) with the cosine of the angle $\varphi$ between the rays going from the origin through the transformed compositions (see Fig.\ \ref{hellinger}), also referred to as Bhattacharyya coefficient \cite{Bcoeff}:
\begin{equation}
    \mathrm{cos}~\varphi=\sum_{i=1}^D\sqrt{x_iy_i}.
\end{equation}

Let us now come back to Aitchison distance and discuss in which structural aspects it differs from divergences obtained from the Fisher metric. In fact, in data analysis, parametrized classes of distance measures are common \cite{green3}. In the same way as in Eq.\ (\ref{alphadiv}), they are mediated by the Box-Cox transformation, which has the limit
\begin{equation}
    \lim_{\beta\to0}\frac{x^\beta-1}{\beta}=\log(x).\label{BoxCox}
\end{equation}
This has been applied in CoDA to obtain log-ratio analysis as a Correspondence Analysis of power-transformed data \cite{green1}. There, we have the following family of distance measures:
\begin{equation}
    d^2_\beta(\boldsymbol{x},\boldsymbol{y})=\frac{1}{\beta^2}\sum_{i=1}^D\omega_i\left(\mathcal{C}{(\boldsymbol{x}^\beta)_i}-\mathcal{C}{(\boldsymbol{y}^\beta)_i}\right)^2,\label{chisquare}
\end{equation}
where $\omega_i$ are suitable weights. For the case $\beta=1$, this is the (square of the) symmetric $\chi^2$-distance used in Correspondence Analysis, while $\beta=0$ gives\footnote{This can be seen as the high-temperature limit in statistical physics.} Aitchison distance (when $\omega_i=D^2$ for all $i$), see the Appendix for a proof. Although the case $\beta=1/2$ has a direct relationship with it, Hellinger distance does not form part of this family because of the closure operation that makes us stay inside the simplex. Similarly, Aitchison distance cannot be obtained from the alpha divergences of Eq.\ (\ref{alphadiv}). Alpha divergences are included in a more general class of divergences known under the name of $f$-divergence. They have the form
\begin{equation}
    D_f(\boldsymbol{x},\boldsymbol{y})=\sum_{i=1}^Dx_if\left(\frac{y_i}{x_i}\right),
\end{equation}
where $f$ is a convex function. It is a well-established result in information geometry that $f$-divergences are the only decomposable\footnote{A divergence is decomposable if it can be written as a sum of terms that only depend on individual components.} divergences that behave monotonically under coarse-graining of information \cite{Amari}, i.e., when compositional parts are amalgamated into higher-level parts. This important invariance property is called information monotonicity. Aitchison distance is not decomposable, as in each summand we use information from all compositional parts when evaluating their geometric mean. Nevertheless, in the next section we will show that it fulfills information monotonicity.\\
It is interesting to note that the two Euclidean distances (Hellinger and Aitchison) are each related to different isometries of the simplex. In the case of Hellinger distance, compositions are isometrically mapped into the positive orthant of the sphere\footnote{But note that this mapping does not obtain the spherical representative of the composition in the sense of the definition of an equivalence class.}. This isometry also holds for the Fisher metric itself \cite{Nihat} (which makes it possible to view the Fisher metric as a Euclidean metric on the $\mathbb{R}^D$ where the sphere is embedded). In the case of Aitchison distance, the isometry in question is of central interest in CoDA. It is the clr transformation, i.e., the map between $\mathcal{S}^D$ and $\mathcal{T}^D$. Here, however, there is no corresponding isometry of the Fisher metric. Although a Euclidean metric may appear convenient, it does not have the same desirable properties as the Fisher metric. In fact, a central result in information geometry states that the Fisher metric is the only metric \cite{Chentsov2} that stays invariant under coarse graining of information. 

\section{Information monotonicity of Aitchison distance}
\subsection{Amalgamations lead to coarse grained information}\label{coarse}
 Let us denote by $\mathcal{A}$ a subset of $\mathcal{D}:=\{1,\dots,D\}$, and let $\boldsymbol{x}_{\mathcal{A}}$ denote the corresponding subcomposition, i.e., the vector where parts with indices not belonging to $\mathcal{A}$ have been removed. Let $H(\boldsymbol{x})$ denote the Shannon entropy of composition $\boldsymbol{x}$, i.e., the potential $-\phi(\boldsymbol{\eta}_{\boldsymbol{x}})$. Consider its decomposition
\begin{equation}
    H(\boldsymbol{x})=(1-s(\boldsymbol{x}_{\mathcal{A}}))H(\mathcal{C}\boldsymbol{x}_{\mathcal{D}\backslash\mathcal{A}})+s(\boldsymbol{x}_{\mathcal{A}})H(\mathcal{C}\boldsymbol{x}_{\mathcal{A}})+H(s(\boldsymbol{x}_{\mathcal{A}}),1-s(\boldsymbol{x}_{\mathcal{A}})),\label{entropySub}
\end{equation}
where $s(\boldsymbol{x}_{\mathcal{A}})=\sum_{i\in\mathcal{A}}x_i$. We see that this is a convex combination of the entropies of the two subcompositions plus a binary entropy, where all terms involve the amalgamation $s(\boldsymbol{x}_{\mathcal{A}})$. Probabilistically speaking, this particular amalgamation corresponds to the probability that an event occurs for any of the events $\mathcal{A}$ that are left out to obtain the subcomposition $\boldsymbol{x}_{\mathcal{D}\backslash\mathcal{A}}$. Generally, amalgamation is nothing else but a coarse graining of the events and their probabilities. The corresponding coarse graining of information is described by this alternative decomposition of Shannon entropy:
\begin{equation}
    H(\boldsymbol{x})=H(\boldsymbol{x}_{\mathcal{D}\backslash\mathcal{A}},s(\boldsymbol{x}_{\mathcal{A}}))+s(\boldsymbol{x}_{\mathcal{A}})H(\mathcal{C}\boldsymbol{x}_\mathcal{A}).\label{entropyAm}
\end{equation}
As the second summand is greater equal zero, this also shows that information cannot grow after coarse graining.

\subsection{The notion of monotonicity for divergences and distance measures}
These considerations lead us to an important result that concerns the divergence associated with the potential $\phi(\boldsymbol{\eta}_{\boldsymbol{x}})$, i.e., the relative entropy. This result is the information monotonicity under coarse graining, where the notion of monotonicity is somewhat related to the notion of subcompositional dominance. The latter refers to the property that a measure of distance does not increase when evaluating it on a subset of parts only. This is often seen as a desirable property of distances in CoDA (and not fulfilled by distances like Hellinger and Batthacharyya, see \cite{clustering} for a discussion of distance measures with respect to compositions). A similar---but perhaps more natural---requirement that has not received attention yet in the CoDA community is the one that a distance between compositions should not increase when comparing it with the one after amalgamating over a subset of parts.\footnote{Subcompositional coherence, i.e., the fundamental requirement that quantities remain identical on a renormalized subcomposition, is not an issue for amalgamation: after amalgamation there is no need for renormalization.} As we have seen in the previous subsection, we cannot gain information when amalgamating parts, so we should lose resolution when comparing the amalgamated compositions. This is also related to the notion of sufficient statistic, see \cite{Amari}. The information monotonicity property of relative entropy can be expressed as     
\begin{equation}
    D_\phi(\boldsymbol{x}||\boldsymbol{y})\ge D_\phi\left((\boldsymbol{x}_{\mathcal{D}\backslash\mathcal{A}},s(\boldsymbol{x}_{\mathcal{A}}))~||~(\boldsymbol{y}_{\mathcal{D}\backslash\mathcal{A}},s(\boldsymbol{y}_{\mathcal{A}}))\right).\label{infomono}
\end{equation}
This result can be shown for the more general case of $f$-divergences and continuous distributions using Jensen's inequality \cite{AmariNagaoka}.\\
Note that in \cite{balances}, when discussing amalgamations of parts, the notion of ``monotonicity'' is used differently. There, the authors argue against amalgamations referring to the observation that Aitchison distances between amalgamated compositions and the amalgamated center of the simplex show a non-monotonic behaviour along an ilr-coordinate axis defined before amalgamation. We will show below that information monotonicity does hold for Aitchison distance. We see the lack of distance monotonicity as discussed in \cite{balances} rather as an argument against the use of a Euclidean coordinate system here. 

\subsection{Monotonicity of Aitchison distance}
A symmetrized version of relative entropy has recently been used in the context of data-driven amalgamation \cite{amalgams}, where it was shown to be better preserved between samples than Aitchison distance. While the information-theoretic meaning and mathematical properties reflected in the decompositions shown in section \ref{coarse} make Shannon entropy an ideal measure of information, alternative indices that can sometimes be evaluated more easily on real-world data (e.g., making use of sums of squares) have also been considered. In our context, it is interesting that $A(\boldsymbol{z}_{\boldsymbol{x}})$, the potential of Eq.\ (\ref{ilrpot}), has been proposed as an alternative measure of information within an attempt to reformulate information theory from a CoDA point-of-view \cite{E&PGinfo}. More recently, it has also been proposed as an inequality index (when divided by the number of parts $D$) \cite{sampleSpace}. Here, the following decomposition was shown:
\begin{equation}
    ||\boldsymbol{x}||_A^2=||(\boldsymbol{x}_{\mathcal{D}\backslash\mathcal{A}},\boldsymbol{x}_{\mathcal{A}})||_A^2=||\boldsymbol{x}_{\mathcal{D}\backslash\mathcal{A}}||_A^2+||\boldsymbol{x}_{\mathcal{A}}||_A^2+\frac{a(D-a)}{D}\left(\log\frac{g(\boldsymbol{x}_{\mathcal{D}\backslash\mathcal{A}})}{g(\boldsymbol{x}_{\mathcal{A}})}\right)^2,\label{clrpot}
\end{equation}
where we denoted the set size of $\mathcal{A}$ by $a$. If we now want to decompose with respect to a composition that was partly amalgamated, we find a corresponding relationship that is more complicated\footnote{It is interesting to note that the two interaction terms have the form of squares of the balance and pivot coordinates mentioned in section \ref{CoDA}.}:
\begin{multline}
    ||\boldsymbol{x}||_A^2=||(\boldsymbol{x}_{\mathcal{D}\backslash\mathcal{A}},s(\boldsymbol{x}_{\mathcal{A}})||_A^2+||\boldsymbol{x}_{\mathcal{A}}||_A^2\\
    +\frac{a(D-a)}{D}\left(\log\frac{g(\boldsymbol{x}_{\mathcal{D}\backslash\mathcal{A}})}{g(\boldsymbol{x}_{\mathcal{A}})}\right)^2-\frac{D-a}{D-a+1}\left(\log\frac{g(\boldsymbol{x}_{\mathcal{D}\backslash\mathcal{A}})}{s(\boldsymbol{x}_{\mathcal{A}})}\right)^2.
\end{multline}
Clearly, if we replace the amalgamation $s(\boldsymbol{x}_{\mathcal{A}})$ by the geometric mean $g(\boldsymbol{x}_{\mathcal{A}})$, we get a simpler equality:
\begin{multline}
    ||\boldsymbol{x}||_A^2=||(\boldsymbol{x}_{\mathcal{D}\backslash\mathcal{A}},g(\boldsymbol{x}_{\mathcal{A}}))||_A^2+||\boldsymbol{x}_{\mathcal{A}}||_A^2\\
    +\left(\frac{a(D-a)}{D}-\frac{D-a}{D-a+1}\right)\left(\log\frac{g(\boldsymbol{x}_{\mathcal{D}\backslash\mathcal{A}})}{g(\boldsymbol{x}_{\mathcal{A}})}\right)^2.\label{geom}
\end{multline}
 Aggregating by geometric means or by amalgamations has been a subject of debate in the CoDA community \cite{sampleSpace, amalgGreen}. As we can see, measures like the Aitchison norm lend themselves much better to taking geometric means rather than amalgamations. There is, however, no straight-forward probabilistic interpretation of geometric means\footnote{The product over parts specifies the probability that {\it all} events in the subset occur, but this is then re-scaled by the exponent to the probability of a single event.}, and the more elegant formal expressions that result often suffer from reduced interpretability.\\
To the best of our knowledge, the information monotonicity property in its general form has not been considered yet for Aitchison distance. We here exploit the various decompositions stated above for proving it. Results are summarized in the following propositions. All proofs can be found in the Appendix.
{\proposition
Let $\mathcal{D}$ and $\mathcal{A}\subset\mathcal{D}$ be two index sets with sizes $D$ and $a$, respectively. Further, let $\boldsymbol{x}$ and $\boldsymbol{y}$ be the simplicial representatives of two compositions in $\mathcal{S}^D$. The amalgamation of $\boldsymbol{x}$ over the subset $\mathcal{A}$ of parts be denoted by  $s(\boldsymbol{x}_{\mathcal{A}})=\sum_{i\in\mathcal{A}}x_i$. Then the following decomposition of Aitchison distance holds:
\begin{multline}
    ||\boldsymbol{x}\ominus\boldsymbol{y}||_A^2=||(\boldsymbol{x}_{\mathcal{D}\backslash\mathcal{A}},s(\boldsymbol{x}_{\mathcal{A}}))\ominus(\boldsymbol{y}_{\mathcal{D}\backslash\mathcal{A}},s(\boldsymbol{y}_{\mathcal{A}}))||_A^2+||\boldsymbol{x}_{\mathcal{A}}\ominus\boldsymbol{y}_{\mathcal{A}}||_A^2\\
    +\frac{a(D-a)}{D}\left(\log\frac{g(\boldsymbol{x}_{\mathcal{D}\backslash\mathcal{A}})}{g(\boldsymbol{x}_{\mathcal{A}})}-\log\frac{g(\boldsymbol{y}_{\mathcal{D}\backslash\mathcal{A}})}{g(\boldsymbol{y}_{\mathcal{A}})}\right)^2\\
    -\frac{D-a}{D-a+1}\left(\log\frac{g(\boldsymbol{x}_{\mathcal{D}\backslash\mathcal{A}})}{s(\boldsymbol{x}_{\mathcal{A}})}-\log\frac{g(\boldsymbol{y}_{\mathcal{D}\backslash\mathcal{A}})}{s(\boldsymbol{y}_{\mathcal{A}})}\right)^2.\nonumber
\end{multline}
}
{\corollary\label{geomeanDecom}
When aggregating a subset of parts in form of their geometric mean, we have the following decomposition of Aitchison distance
\begin{multline}
||\boldsymbol{x}\ominus\boldsymbol{y}||_A^2=||(\boldsymbol{x}_{\mathcal{D}\backslash\mathcal{A}},g(\boldsymbol{x}_{\mathcal{A}}))\ominus(\boldsymbol{x}_{\mathcal{D}\backslash\mathcal{A}},g(\boldsymbol{y}_{\mathcal{A}}))||_A^2+||\boldsymbol{x}_{\mathcal{A}}\ominus\boldsymbol{y}_{\mathcal{A}}||_A^2\\
    +\left(\frac{a(D-a)}{D}-\frac{D-a}{D-a+1}\right)\left(\log\frac{g(\boldsymbol{x}_{\mathcal{D}\backslash\mathcal{A}})}{g(\boldsymbol{x}_{\mathcal{A}})}-\log\frac{g(\boldsymbol{y}_{\mathcal{D}\backslash\mathcal{A}})}{g(\boldsymbol{y}_{\mathcal{A}})}\right)^2.
\end{multline}
}

From this decomposition, we get the following monotonicity result: 
{\corollary\label{geomeanIneq}
With parts aggregated by geometric means, the following inequality holds
\[||\boldsymbol{x}\ominus\boldsymbol{y}||_A^2\ge||(\boldsymbol{x}_{\mathcal{D}\backslash\mathcal{A}},g(\boldsymbol{x}_{\mathcal{A}}))\ominus(\boldsymbol{x}_{\mathcal{D}\backslash\mathcal{A}},g(\boldsymbol{y}_{\mathcal{A}}))||_A^2+||\boldsymbol{x}_{\mathcal{A}}\ominus\boldsymbol{y}_{\mathcal{A}}||_A^2.\]
}
As we can see, for geometric-mean summaries, the sum of the interaction terms (i.e.,  of the terms not involving norms) will remain greater equal zero. This is no longer true for amalgamation of parts, and it is less straightforward to show the corresponding inequality:
{\proposition\label{infomono}
Aitchison distance fulfills the information monotonicity
\[||\boldsymbol{x}\ominus\boldsymbol{y}||_A^2\ge||(\boldsymbol{x}_{\mathcal{D}\backslash\mathcal{A}},s(\boldsymbol{x}_{\mathcal{A}}))\ominus(\boldsymbol{x}_{\mathcal{D}\backslash\mathcal{A}},s(\boldsymbol{y}_{\mathcal{A}}))||_A^2.\]
}
\\

\section{Discussion and outlook}
In our little outline of finite information geometry, we could but scratch the surface of the formal apparatus that is at our disposal. We are certain it can serve to advance the field of Compositional Data Analysis in various ways. Differential geometry provides a universally valid framework for the problems occurring with constrained data. Considering the simplex a differentiable manifold enables a general approach from which specific problems like the compositional differential calculus \cite{compDiffCalc} follow naturally. Clearly, there is (and has to be) overlap in methodology between the information-geometric perspective and the CoDA approach. An example is the use of the exponential map anchored at the center of the simplex discussed in section \ref{divergences}, which allows to identify the simplex with a linear space that is central to the Euclidean CoDA approach. Another example is the fundamental role that exponential families play in information geometry and which have been studied in the CoDA context in the so-called Bayes spaces \cite{BayesSpaces}. But we also think that some of the current limitations of CoDA can be overcome using the additional structures that information geometry can provide. The ease with which amalgamations of parts can be handled by Kullback-Leibler divergence might partly resolve the debate surrounding this issue in the CoDA community. Further, maximum-entropy projections, where Kullback-Leibler divergences are the central tool, seem an especially promising avenue to pursue in the context of data that are only partially available or subject to constraints.  Also, our description has focused on the equivalence of compositions with discrete probability distributions, but information geometry can of course be used to describe the distributions themselves. These are no longer finite but continuous and contain a constraint that introduces dependencies among their random variables, calling for the use of more general versions of the concepts presented here. 

\section*{Appendix}
\subsubsection*{Proof that Eq.\ \ref{chisquare} tends to Aitchison distance for $\beta\rightarrow 0$}
\begin{equation}
    d^2_\beta(\boldsymbol{x},\boldsymbol{y})=\frac{1}{\beta^2}\sum_{i=1}^D\omega_i\left(\mathcal{C}{(\boldsymbol{x}^\beta)_i}-\mathcal{C}{(\boldsymbol{y}^\beta)_i}\right)^2=\sum_{i=1}^D\omega_i\left(\frac{x^\beta_i}{\beta\sum_{k=1}^Dx^\beta_k}-\frac{y^\beta_i}{\beta\sum_{k=1}^Dy^\beta_k}\right)^2.\nonumber
\end{equation}
The terms inside the bracket can be written as
\begin{equation}
    \frac{x^\beta_i}{\beta\sum_{k=1}^Dx^\beta_k}-\frac{y^\beta_i}{\beta\sum_{k=1}^Dy^\beta_k}=\frac{x^\beta_i-\frac{1}{D}\sum_{k=1}^Dx^\beta_k}{\beta\sum_{k=1}^Dx^\beta_k}-\frac{y^\beta_i-\frac{1}{D}\sum_{k=1}^Dy^\beta_k}{\beta\sum_{k=1}^Dy^\beta_k}
\end{equation}
when subtracting $1/(\beta D)$ and adding it back. Similarly,
\begin{equation}
    =\frac{x^\beta_i-1-\frac{1}{D}\sum_{k=1}^D(x^\beta_k-1)}{\beta\sum_{k=1}^Dx^\beta_k}-\frac{y^\beta_i-1-\frac{1}{D}\sum_{k=1}^D(y^\beta_k-1)}{\beta\sum_{k=1}^Dy^\beta_k}.
\end{equation}
In this, we recognize the Box-Cox transform in the numerators, and can make use of the limit in Eq.\ (\ref{BoxCox}). The sums in the denominators clearly tend to $D$ for $\beta\rightarrow 0$. We can thus evaluate the limit as a quotient of finite limits and conclude
\begin{equation}
    \lim_{\beta\to0}d^2_\beta(\boldsymbol{x},\boldsymbol{y})=\sum_{i=1}^D\frac{\omega_i}{D^2}\left(\log\frac{x_i}{g(\boldsymbol{x})}-\log\frac{y_i}{g(\boldsymbol{y})}\right)^2.
\end{equation}

\subsubsection*{Proof of Proposition 1}

We start by defining $\boldsymbol{z}=(x_i/y_i)_{i\in\mathcal{D}}$. We can now use the decomposition 
\begin{equation}
||\boldsymbol{z}||_A^2=||\boldsymbol{z}_{\mathcal{D}\backslash\mathcal{A}}||_A^2+||\boldsymbol{z}_{\mathcal{A}}||_A^2+\frac{a(D-a)}{D}\left(\log\frac{g(\boldsymbol{z}_{\mathcal{D}\backslash\mathcal{A}})}{g(\boldsymbol{z}_{\mathcal{A}})}\right)^2,\label{decom}
\end{equation}
which can be derived using equalities like
\begin{equation}
    g(\boldsymbol{x})=\left(g(\boldsymbol{x}_{\mathcal{A}})^ag(\boldsymbol{x}_{\mathcal{D}\backslash\mathcal{A}})^{D-a}\right)^{\frac{1}{D}}=g(\boldsymbol{x}_{\mathcal{A}})g(\boldsymbol{x}_{\mathcal{A}})^{\frac{a}{D}-1}g(\boldsymbol{x}_{\mathcal{D}\backslash\mathcal{A}})^{\frac{D-a}{D}},\nonumber
\end{equation} 
which can be used to expand
\begin{equation}
    \sum_{i\in\mathcal{A}}
    \left(\log\frac{z_i}{g(\boldsymbol{z})}\right)^2
    =\sum_{i\in\mathcal{A}}\left(\log\frac{z_i}{g(\boldsymbol{z}_{\mathcal{A}})}+\frac{D-a}{D}\log\frac{g(\boldsymbol{z}_{\mathcal{A}})}{g(\boldsymbol{z}_{\mathcal{D}\backslash\mathcal{A}})}\right)^2\nonumber
\end{equation}
and doing the square, cross terms vanish after summation.\\
Next, we observe that, for an arbitrary $s_{\boldsymbol{z}}$ which we join as an additional component with the vector $\boldsymbol{z}_{\mathcal{D}\backslash\mathcal{A}}$, we have
\begin{equation}
 ||(\boldsymbol{z}_{\mathcal{D}\backslash\mathcal{A}},s_{\boldsymbol{z}})||_A^2  = ||\boldsymbol{z}_{\mathcal{D}\backslash\mathcal{A}}||_A^2+\frac{D-a}{D-a+1}\left(\log\frac{g(\boldsymbol{z}_{\mathcal{D}\backslash\mathcal{A}})}{s_{\boldsymbol{z}}}\right)^2\label{amalgnorm}
\end{equation}
because, similarly as before, we have
\begin{multline}
    \sum_{i\in\mathcal{D}\backslash\mathcal{A}}\left(\log\frac{z_i}{g((\boldsymbol{z}_{\mathcal{D}\backslash\mathcal{A}},s_{\boldsymbol{z}}))}\right)^2+\left(\log\frac{s_{\boldsymbol{z}}}{g((\boldsymbol{z}_{\mathcal{D}\backslash\mathcal{A}},s_{\boldsymbol{z}}))}\right)^2=\\
    \sum_{i\in\mathcal{D}\backslash\mathcal{A}}\left(\log\frac{z_i}{g(\boldsymbol{z}_{\mathcal{D}\backslash\mathcal{A}})}+\frac{1}{D-a+1}\log\frac{g(\boldsymbol{z}_{\mathcal{D}\backslash\mathcal{A}})}{s_{\boldsymbol{z}}}\right)^2+\left(
    \frac{D-a}{D-a+1}\log\frac{s_{\boldsymbol{z}}}{g(\boldsymbol{z}_{\mathcal{D}\backslash\mathcal{A}})}\right)^2.\nonumber
\end{multline}

We can now choose for $s_z$ the expression $s_z=(\sum_{i\in\mathcal{A}}x_i)/(\sum_{i\in\mathcal{A}}y_i)$. With this, we express
$||\boldsymbol{z}_{\mathcal{D}\backslash\mathcal{A}}||_A^2$ using the expression following from Eq.\ (\ref{amalgnorm}). We then substitute the corresponding term in Eq.\ (\ref{decom}), proving the proposition:
\begin{multline}
    ||\boldsymbol{z}||_A^2=||(\boldsymbol{z}_{\mathcal{D}\backslash\mathcal{A}},s)||_A^2+||\boldsymbol{z}_{\mathcal{A}}||_A^2\\
    +\frac{a(D-a)}{D}\left(\log\frac{g(\boldsymbol{z}_{\mathcal{D}\backslash\mathcal{A}})}{g(\boldsymbol{z}_{\mathcal{A}})}\right)^2-\frac{D-a}{D-a+1}\left(\log\frac{g(\boldsymbol{z}_{\mathcal{D}\backslash\mathcal{A}})}{(\sum_{i\in\mathcal{A}}x_i)/(\sum_{i\in\mathcal{A}}y_i)}\right)^2.\label{propo}
\end{multline}

\subsubsection*{Proof of Corollary \ref{geomeanDecom}}
This follows from the fact that we can insert $g(\boldsymbol{z}_{\mathcal{A}})$ for $s_{\boldsymbol{z}}$ in Eq.\ (\ref{amalgnorm}). This gets us an expression for $||\boldsymbol{z}_{\mathcal{D}\backslash\mathcal{A}}||_A^2$, which is inserted in Eq.\ (\ref{decom}).

\subsubsection*{Proof of Corollary \ref{geomeanIneq}}
To prove the corollary, we have to show that the last term in the decomposition of Corollary \ref{geomeanDecom} is greater or equal to zero. We thus need to show 
\begin{equation}
    (D-a)\frac{a(D-a+1)-D}{D(D-a+1)}\ge0.\label{FofA}
\end{equation}
Since $D>a>1$, and the quadratic equation $a^2-a(D+1)+D=0$ has solutions $a=1$ and $a=D$, between these values we are either above or below zero. We are above because the first derivative of Eq.\ (\ref{FofA}) in $a=1$ is $(D-1)^2/D^2$, which is bigger zero.\\

\subsubsection*{Proof of Proposition \ref{infomono}}
Let $\boldsymbol{z}$ denote the vector with components $x_i/y_i$, $i=1,\dots,D$. To prove the proposition, we have to bound the terms after the first plus sign in Proposition 1 from below, i.e.,
\begin{multline}
    R(\boldsymbol{x},\boldsymbol{y}):=\sum_{i\in\mathcal{A}}\left(\log\frac{z_i}{g(\boldsymbol{z}_{\mathcal{A}})}\right)^2+\frac{a(D-a)}{D}\left(\log\frac{g(\boldsymbol{z}_{\mathcal{D}\backslash\mathcal{A}})}{g(\boldsymbol{z}_{\mathcal{A}})}\right)^2\\
    -\frac{D-a}{D-a+1}\left(\log\frac{g(\boldsymbol{z}_{\mathcal{D}\backslash\mathcal{A}})}{(\sum_{i\in\mathcal{A}}x_i)/(\sum_{i\in\mathcal{A}}y_i)}\right)^2\ge0.
\end{multline}
Let us start with the second summand. We rewrite it as
\begin{equation}
    \left(\frac{a(D-a)}{D}-\frac{D-a}{D-a+1}\right)\left(\log\frac{g(\boldsymbol{z}_{\mathcal{D}\backslash\mathcal{A}})}{g(\boldsymbol{z}_{\mathcal{A}})}\right)^2+\frac{D-a}{D-a+1}\left(\log\frac{g(\boldsymbol{z}_{\mathcal{D}\backslash\mathcal{A}})}{g(\boldsymbol{z}_{\mathcal{A}})}\right)^2\nonumber
\end{equation}
Since we showed Eq.\ (\ref{FofA}), the summand on the left is greater equal zero. Thus we have
\begin{multline}
    R(\boldsymbol{x},\boldsymbol{y})\ge\sum_{i\in\mathcal{A}}\left(\log\frac{z_i}{g(\boldsymbol{z}_{\mathcal{A}})}\right)^2\\
    +\frac{D-a}{D-a+1}\left(\left(\log\frac{g(\boldsymbol{z}_{\mathcal{D}\backslash\mathcal{A}})}{g(\boldsymbol{z}_{\mathcal{A}})}\right)^2-\left(\log\frac{g(\boldsymbol{z}_{\mathcal{D}\backslash\mathcal{A}})}{(\sum_{i\in\mathcal{A}}x_i)/(\sum_{i\in\mathcal{A}}y_i)}\right)^2\right)\\
    \ge\sum_{i\in\mathcal{A}}\left(\log\frac{z_i}{g(\boldsymbol{z}_{\mathcal{A}})}\right)^2-\left(\log\frac{(\sum_{i\in\mathcal{A}}x_i)/(\sum_{i\in\mathcal{A}}y_i)}{g(\boldsymbol{z}_{\mathcal{A}})}\right)^2,\label{rest}
\end{multline}
where the second inequality follows because the big bracket (with a prefactor smaller one) has a structure that can be bounded like
\begin{equation}
    |(A-B)^2-(A-C)^2|\le(B-C)^2,\nonumber
\end{equation}
with $\log~g(\boldsymbol{z}_{\mathcal{D}\backslash\mathcal{A}})$ playing the role of $A$. Finally, the last term in Eq.\ (\ref{rest}) can be bounded from above as follows. Since $x_i\le y_i\cdot\max x_i/y_i$, we also have
\begin{equation}
    \sum_{i\in\mathcal{A}}x_i\le\max_{k\in\mathcal{A}}\frac{x_k}{y_k}\sum_{i\in\mathcal{A}}y_i,\label{maxsum}
\end{equation}
so the ratio of sums is smaller than the maximum over the ratios. Without restricting generality, let us assume that $g(\boldsymbol{x}_{\mathcal{A}})\ge g(\boldsymbol{y}_{\mathcal{A}})$, i.e., $g(\boldsymbol{z}_{\mathcal{A}})\ge1$. The bound on the sum ratio implied by Eq.\ (\ref{maxsum}) is then sufficient for proving the proposition:
\begin{equation}
    R(\boldsymbol{x},\boldsymbol{y})\ge\max_{i\in\mathcal{A}}\left(\log\frac{x_i/y_i}{g(\boldsymbol{z}_{\mathcal{A}})}\right)^2-\left(\log\frac{(\sum_{i\in\mathcal{A}}x_i)/(\sum_{i\in\mathcal{A}}y_i)}{g(\boldsymbol{z}_{\mathcal{A}})}\right)^2\ge0.
\end{equation}


\begin{thebibliography}{99.}%

\bibitem{AitchisonBook} Aitchison, J: The statistical analysis of compositional data. Chapman and Hall (1986)

\bibitem{green} Greenacre, M: Compositional Data Analysis, Annual Reviews of Statistics and its Application \textbf{8}(1), 271--299 (2021)

\bibitem{sampleSpace} Egozcue, JJ and Pawlowsky-Glahn, V: Compositional data: the sample space and its structure. TEST \textbf{28}(3), 599--638 (2019)

\bibitem{mathCoDA} Barcel{\'o}-Vidal, C, Mart{\'\i}n-Fern{\'a}ndez, JA: The Mathematics of Compositional Analysis. Austrian Journal of Statistics \textbf{45}(4), 57--71 (2016)

\bibitem{logratioTrans} Aitchison, J: The Statistical Analysis of Compositional Data. J Royal Stat Soc B \textbf{44} (2), 139--160 (1982)

\bibitem{simplicialGeometry} Egozcue, JJ, Barcel{\'o}-Vidal, C, Mart{\'\i}n-Fern{\'a}ndez, JA, Jarauta-Bragulat, E, D{\'\i}az-Barrero, JL, Mateu-Figueras, G, Pawlowsky-Glahn, V, Buccianti, A: Elements of simplicial linear algebra and geometry. In: Pawlowsky-Glahn, V and Buccianti, A (eds.) Compositional data analysis: Theory and applications, pp. 141--157. Wiley (2011)

\bibitem{ilr} Egozcue, JJ, Pawlowsky-Glahn, V, Mateu-Figueras, G,
and Barcel\'o-Vidal, C: Isometric Logratio Transformations for Compositional Data Analysis. Mathematical Geology, \textbf{35}  (3), 279--300 (2003)

\bibitem{balances} J. J. Egozcue and V. Pawlowsky-Glahn. Groups of Parts and Their Balances in Compositional Data Analysis. Mathematical Geology, \textbf{37} (7), 795--828 (2005)

\bibitem{pivot} Hron, K, Filzmoser, P, de Caritat, P, Fi\^serov\'a, E, Gardlo, A: Weighted Pivot Coordinates for Compositional Data and Their Application to Geochemical Mapping. Mathematical Geosciences \textbf{49}, 797--814 (2017) 

\bibitem{Chentsov} Chentsov, N: Statistical Decision Rules and Optimal Inference (vol.\ 53), Nauka (1972) (in Russian); English translation in: Math. Monograph. (vol.\ 53), Am. Math. Soc. (1982)

\bibitem{Rao} Rao, CR: Information and the accuracy attainable in the estimation of statistical parameters.
Bull. Calcutta Math. Soc. \textbf{37}, 81--89 (1945)

\bibitem{AmariNagaoka} Amari, S and Nagaoka, H: Methods of Information Geometry. Translations of Mathematical Monographs (vol.\ 191), American Mathematical Society (2000) 

\bibitem{lectureNotes} Amari, S: Differential-Geometric Methods in Statistics. Lecture Notes in Statistics (vol.\ 28),
Springer (1985)

\bibitem{Nihat} Ay, N, Jost, J, Le, HV, Schwachhöfer, L: Information Geometry. A Series of Modern Surveys
in Mathematics (vol.\ 64), Springer (2017)

\bibitem{Amari} Amari, S: Information Geometry and Its Applications. Applied Mathematical Sciences,
(vol.\ 194), Springer (2016)

\bibitem{graphMod} Whittaker, J: Graphical models in applied multivariate statistics. Wiley (1990)

\bibitem{AyErb} Ay, N and Erb, I: On a notion of linear replicator equations. J Dyn. Diff. Eqs., \textbf{17} (2), 427-451 (2005)

\bibitem{symmDiv} Mart\'in-Fern\'andez, JA, Bren, M, Barcel\'o-Vidal, C and  Pawlowsky-Glahn, V: A Measure of Difference for Compositional Data based on measures of divergence. In: Proceedings of the Fifth Annual  Conference of the International Assotiation for Mathematical Geology, Ed.  Lippard, S.J., Naess, A., and  Sinding-Larsen, R.. Trondheim (Norway), Vol. 1, pp. 211-215 (1999)

\bibitem{Renyi} R\'enyi, A.: On measures of entropy and information. In: Proceedings of the 4th Berkeley
Symposium on Mathematical Statistics and Probability, vol. 1, pp.\ 547–561. University of California Press, Berkeley (1961)

\bibitem{Tsallis} Tsallis, C.: Possible generalization of Boltzmann–Gibbs statistics. J. Stat. Phys. \textbf{52}, 479–487
(1988)

\bibitem{green3} Greenacre, M: Power transformations in correspondence analysis. Computational Statistics \& Data Analysis, \textbf{53}(8), 3107-3116 (2009)

\bibitem{green1} Greenacre, M: Log-Ratio Analysis Is a Limiting Case of Correspondence Analysis. Math Geosci \textbf{42}, 129 (2010) 

\bibitem{Bcoeff} Bhattacharyya, A: On a measure of divergence between two statistical populations defined by their probability distributions. Bulletin of the Calcutta Mathematical Society \textbf{35}, 99–109. (1943)

\bibitem{Chentsov2} Chentsov, N: Algebraic foundation of mathematical statistics. Math.\ Oper.forsch.\ Stat., Ser.\ Stat.\ \textbf{9}, 267–276 (1978)

\bibitem{clustering} Palarea-Albaladejo, J, Martín-Fern\'andez, JA and Soto, J.A.:  Dealing with Distances and Transformations for Fuzzy C-Means Clustering of Compositional Data, Journal of Classification \textbf{29}(2), 144-169 (2012)

\bibitem{amalgams} Quinn, TP, Erb, I: Amalgams: data-driven amalgamation for the dimensionality reduction of compositional data, NAR Genomics and Bioinformatics \textbf{2}(4) Iqaa076 (2021)

\bibitem{E&PGinfo} Egozcue, JJ and Pawlowsky-Glahn, V: Evidence functions: a compositional approach to information. SORT \textbf{42} (2), 101-124 (2018)

\bibitem{amalgGreen} Greenacre, M: Amalgamations are valid in compositional data analysis, can be used in agglomerative clustering, and their logratios have an inverse transformation. Appl.\ Comp.\ Geosc., \textbf{5}, 100017 (2020)

\bibitem{compDiffCalc} Barcel\'o-Vidal, C., Mart\'in-Fern\'andez, J.A. and Mateu-Figueras, G.,  Compositional Differential Calculus on the Simplex. In: Pawlowsky-Glahn, V and Buccianti, A (eds.) Compositional data analysis: Theory and applications, pp. 176--190. Wiley (2011)

\bibitem{BayesSpaces} Egozcue, JJ,  Pawlowsky-Glahn, V,  Tolosana-Delgado, R, Ortego, MI  and van  den  Boogaart, KG,  Bayes  spaces:  use  of  improper  distributions  and exponential families. Revista de la Real Academia de Ciencias Exactas, F\'isicas y Naturales. Serie A. Matematicas, \textbf{107}(2), pp. 475--486 (2013)


%
\bigskip

\end{thebibliography}
\end{document}